\documentclass[a4paper,12pt]{amsart}
\theoremstyle{definition}
\newtheorem{definition}{Definition}[section]

\theoremstyle{plain}
\newtheorem{theorem}[definition]{Theorem}
\newtheorem{corollary}[definition]{Corollary}
\newtheorem{lemma}[definition]{Lemma}
\newtheorem{proposition}[definition]{Proposition}

\newtheorem*{theorem*}{Theorem}
\newtheorem*{lemma*}{Lemma}
\newtheorem*{proposition*}{Proposition}
\newtheorem*{problem*}{Problem}
\newtheorem*{conjecture}{Conjecture}
\newtheorem*{definition*}{Definition}
\newtheorem*{corollary*}{Corollary}

\theoremstyle{definition}

\newcommand{\dR}{\mathbb R}
\newcommand{\dP}{\mathbb P}

\newcommand{\cG}{\mathcal G}
\newcommand{\cL}{\mathcal L}

\newcommand{\cH}{\mathcal H}
\newcommand{\cA}{\mathcal A}

\newcommand{\cB}{\mathcal B}

\DeclareMathOperator{\Jac}{Jac}
\DeclareMathOperator{\ex}{ex}
\DeclareMathOperator{\inte}{int}

\DeclareMathOperator{\vol}{Vol}
\DeclareMathOperator{\ar}{Area}
\DeclareMathOperator{\ent}{Ent}

\title{Volume entropy of Hilbert Geometries}
\author{G. Berck, A. Bernig \& C. Vernicos}
\address{D\'epartement de Math\'ematiques, Chemin du mus\'ee 23, 1700 Fribourg, Switzerland}
\email{gautier.berck@unifr.ch}
\email{andreas.bernig@unifr.ch}
\address{Department of Mathematics, Logic House, South Campus, NUI Maynooth, Co. Kildare, Ireland}
\email{Constantin.Vernicos@maths.nuim.ie}
\thanks{
\\ The first two authors were supported
  by the Schweizerischer Nationalfonds grants SNF PP002-114715/1 and 200020-121506 /1.}
\begin{document}
\begin{abstract}
It is shown that the volume entropy of a Hilbert geometry associated to an $n$-dimensional convex body of class $C^{1,1}$ equals $n-1$. To achieve this result, a new projective invariant of convex bodies, similar to the centro-affine area, is constructed. In the case $n=2$, and without any assumption on the boundary, it is shown that the entropy is bounded above by $\frac{2}{3-d} \leq 1$, where $d$ is the Minkowski dimension of the extremal set of $K$. An example of a plane Hilbert geometry with entropy strictly between $0$ and $1$ is constructed.  
\end{abstract}

\maketitle

\section{Introduction}

In his famous 4-th problem, Hilbert asked to characterize metric geometries whose geodesics are straight lines. He constructed a special class of examples, nowadays called \emph{Hilbert geometries} \cite{hilbert1895,hilbert1899}. These geometries have attracted a lot of interest, see for example 
the works of Y.~Nasu~\cite{nasu61}, P.~de~la~Harpe~\cite{delaHarpe93}, A.~Karlsson \& G.~Noskov~\cite{karlsson_guennadi02}, E.~Socie-Methou~\cite{sociemethou04}, T.~Foertsch \& A.~Karlsson~\cite{foertsch_karlsson05}, Y.~Benoist~\cite{benoist06}, B.~Colbois \& C.~Vernicos~\cite{colbois_vernicos07}
and the two complementary surveys by Y.~Benoist~\cite{yvessurvey} and the last named author \cite{vernicos05}.

A Hilbert geometry is a particularly simple metric space on the interior of a 
compact convex set $K$ (see definition below). 
This metric happens to be a complete Finsler metric 
whose set of geodesics contains the straight lines. 
Since the definition of the Hilbert geometry only uses cross-ratios, 
the Hilbert metric is a projective invariant. 
In the particular case where $K$ is an ellipsoid, 
the Hilbert geometry is isometric to the usual hyperbolic space. 

An important part of the above mentioned works, and of older ones, is to study how different or close
to the hyperbolic geometry these geometries can be. For instance, 
if $K$ is not an ellipsoid, then the metric is never Riemannian, see D.C.~Kay~\cite[Corollary 1]{kay67}. 
This last result is actually related to the fact that among all 
finite dimensional normed vector spaces, many notions of curvatures
are only satisfied by the Euclidean spaces (see also P.~Kelly \& L.~Paige~\cite{kp}, 
P.~Kelly \& E.~Strauss~\cite{ks,ks2}). However, 
if $\partial K$ is sufficiently smooth then the flag curvature, an analog of the sectional curvature,  
of the Hilbert metric is constant 
and equals $-1$ , see for example Z.~Shen~\cite[Example 9.2.2]{shen01}. 
Hence a question one can ask is whether or not these geometries behave like negatively
curved Riemanniann manifold. The example of the triangle geometry which is isometric
to a two dimensional normed vector space (see P.~De la~Harpe~ \cite{delaHarpe93}) shows that things are a little more involved
(see also theorems cited below). The present work is partially inspired by the feeling that Hilbert geometries might be thought as geometries with Ricci curvature bounded from below, and focuses on the volume growth of balls.

Unlike the Riemannian case, where there is only one natural choice of volume, there are several good choices of volume on a Finsler manifold. We postpone this issue to section \ref{sec_preliminaries} and fix just one {\it volume} (like the $n$-dimensional Hausdorff measure) for the moment. 

Let $B(o,r)$ be the metric ball of radius $r$ centered at $o$. The volume entropy of $K$ is defined by the following limit (provided it exists)
\begin{equation}\label{defentropy}
 \ent K:=\lim_{r \to \infty} \frac{\log \vol B(o,r)}{r}.
\end{equation}

The entropy does not depend on the particular choice of the
base point $o \in \inte K$ nor on the particular choice of the volume. 
If $h=\ent K$, then $\vol B(o,r)$ behaves roughly as $e^{hr}$. 

It is well-known and easy to prove (see, e.g., S.~Gallot, D.~Hulin \& J.~Lafontaine~\cite[Section III.H]{gallot_sylvestre_hulin}) that the volume of  a ball of radius $r$ in the $n$-dimensional  hyperbolic space is given, with $\omega_n$ the volume of the Euclidean unit ball of dimension $n$,  
by 
\begin{displaymath}
 n \omega_n \int_0^r (\sinh s)^{n-1} ds =O(e^{(n-1)r}).
\end{displaymath}
It follows that the entropy of an ellipsoid equals $n-1$.

In general, it is not known whether the above limit exists. 
If the convex set $K$ is divisible, which means that a discrete subgroup of 
the group of isometries of the Hilbert geometry acts cocompactly, 
then the entropy is known to exist, see Y.~Benoist~\cite{benoist04}. 
If the convex set is sufficiently smooth, e.g., $C^2$ with positive curvature suffices, 
then the entropy exists 
and equals $n-1$ (see the theorem of B.~Colbois \& P.~Verovic below). 
In general, one may define lower 
and upper entropies $\underline{\ent}$, $\overline{\ent}$ 
by replacing the limit in the definition (\ref{defentropy}) 
by $\liminf$ or $\limsup$. 

There is a well-known conjecture 
(whose origin seems difficult to locate) 
saying that the hyperbolic space 
has maximal entropy among all Hilbert geometries of the same dimension. 

\begin{conjecture} 
 For any $n$-dimensional Hilbert geometry, 
\begin{displaymath}
 \overline{\ent} K \leq n-1.
\end{displaymath}
\end{conjecture}

Notice that such a result is a consequence of Bishop's volume comparison
theorem for a complete Riemannian manifold of Ricci curvature bounded by $-(n-1)$ 
(see \cite[theorem 3.101, i)]{gallot_sylvestre_hulin}).

Several particular cases of the conjecture were treated in the literature. The following
one shows that the volume entropy does not characterize the hyperbolic geometry
among all Hilbert geometries.
\begin{theorem*}{(\bf B.~Colbois \& P.~Verovic \cite{colbois_verovic04})}\\
 If $K$ is $C^2$-smooth with strictly positive curvature, 
then the Hilbert metric of $K$ is bi-Lipschitz to the 
hyperbolic metric and therefore 
\begin{displaymath}
\ent K=n-1. 
\end{displaymath}
\end{theorem*}

\begin{theorem*}{(\bf B.~Colbois, C.~Vernicos \& P.~Verovic \cite{colbois_vernicos_verovic08})}\\
The Hilbert metric associated to a plane convex polygone is bi-Lipschitz 
to the Euclidean plane. In particular, its entropy is $0$. 
\end{theorem*}

Instead of taking the volume of balls, another natural choice is to study 
the volume growth of the metric spheres $S(o,r)$. 
One may define a (spherical) entropy by
\begin{equation}\label{defsentropy}
 \ent^s K:=\lim_{r \to \infty} \frac{\log \vol S(o,r)}{r},
\end{equation}
provided the limit exists. In general, one may define upper and lower spherical entropies 
$\overline{\ent}^s K$ and $\underline{\ent}^s K$ 
by replacing the limits in the definition (\ref{defsentropy}) by a $\limsup$ or $\liminf$. 

The following theorem is a spherical version of the theorem of 
B.~Colbois \& P.~Verovic. 
\begin{theorem*}
 {\bf (A.A.~Borisenko \& E.A.~Olin \cite{borisenko_olin07})}\\
If $K$ is an $n$-dimensional convex body of class $C^3$ with positive Gauss curvature, then $\ent^s=n-1$. 
\end{theorem*}

Our first main theorem weakens in a substantial way the assumptions 
in the theorem of B.~Colbois \& P.~Verovic 
and strengthens its conclusions for
not only does it give the 
precise value of the entropy but also the \emph{entropy coefficient}. 
In order to state it, we introduce a projective invariant of convex bodies 
interesting in itself. 

Let $V$ be an $n$-dimensional vector space with origin $o$. Given a convex body $K$ containing $o$ in the interior, we define a positive function $a$ on the boundary by the condition that for $p \in \partial K$ we have $-a(p)p \in \partial K$. The letter $a$ stands for \emph{antipodal}. If $V$ is endowed with a Euclidean scalar product, we let $k(p)$ be the Gauss curvature and $n(p)$ be the outer normal vector at a boundary point $p$ (whenever they are well-defined, which is almost everywhere the case following A.D.~Alexandroff~\cite{alexandrov}). 

\begin{definition*} \label{def_cp_area}
The \emph{centro-projective} area of $K$ is
\begin{equation}
 \cA_p(K):=\int_{\partial K}\frac{\sqrt{k}}{\langle n,p\rangle^{\frac{n-1}{2}}}\left(\frac{2a}{1+a}\right)^{\frac{n-1}{2}}\ dA.
\end{equation}
\end{definition*}

It is not quite obvious (but true, as we shall see) that this definition 
does not depend on the choice of the scalar product. In fact, 
the centro-projective area is invariant under 
\emph{projective transformations} fixing the origin. 
The reader familiar with the theory of valuations 
may notice the similarity with the centro-affine surface area, 
whose definition is the same except that the second factor 
(containing the function $a$) does not appear. We refer to the books by Laugwitz \cite{laugwitz} and Leichtweiss \cite{leichtweiss_book} for more information on affine and centro-affine differential geometry. 

\begin{theorem*}[\bf First Main Theorem]
If $K$ is $C^{1,1}$, then 
\begin{equation} \label{eq_entropy_coefficient_intro}
  \lim_{r \to \infty} \frac{\vol B(o,r)}{\sinh^{n-1}r}=\frac{1}{n-1} \cA_p(K)\neq 0 
\end{equation}
and $\ent K=n-1$. Moreover, without any assumption on $K$, if $\cA_p(K) \neq 0$ then $\underline{\ent} K \geq n-1$. 
\end{theorem*} 

In the two-dimensional case the $C^{1,1}$ assumption is not required, indeed, 
we are able to give an upper bound of the
entropy depending on the Minkowski dimension of the set  $\ex K$ of extremal points of $K$.
Recall that an extremal point of a convex body $K$ is a 
point which can not be written as $\frac{a+b}{2}$ with $a,b \in K, a \neq b$. 

\begin{theorem*}[\bf Second Main Theorem] 
Let $K$ be a two-dimensional convex body. Let $d$ be the upper Minkowski dimension of the set
of extremal points of $K$. Then the entropy of $K$ is bounded by 
\begin{equation}
 \overline{\ent K} \leq \frac{2}{3-d} \leq 1.
\end{equation}
Moreover, the equality in \eqref{eq_entropy_coefficient_intro} 
holds true (with $n=2$). 
\end{theorem*}

The inequality is sharp if $K$ is smooth or contains some positively curved smooth part in 
the boundary. In this case the upper Minkowski dimension of $\ex K$ and the entropy both are $1$. On the other hand, for polygones the upper Min\-kow\-ski dimension of the set of extremal points and the entropy both vanish 
(see the theorem of B.~Colbois, C.~Vernicos \& P.~Verovic above), and the inequality is not sharp in this case.  

It should be noted that the entropy behaves in a rather subtle way (see also C.~Vernicos \cite{vernicos08} for a technical and complementary study, to this paper, of the entropy). As we have seen above, the entropy of a polygon vanishes. In contrast to this, we will construct a convex body with piecewise affine boundary whose entropy is strictly between $0$ and~$1$.  
 
Our next theorem, together with the previous ones, shows in particular that it suffices to assume $K$ to be merely of class $C^{1,1}$ in the theorem of A.A.~Borisenko \& E.A.~Olin. 
\begin{theorem*}
For each convex body $K$, 
\begin{align*}
 \underline{\ent}^s K & =\underline{\ent} K,\\
 \overline{\ent}^s K &=\overline{\ent} K. 
\end{align*}
\end{theorem*}

\subsection*{Plan of the paper}
In the next section, we collect some well-known facts about convex bodies, Hilbert geometries and volumes on Finsler manifolds. A number of easy lemmas is proved which will be needed in the proof of our main theorem. Using some inequalities for volumes in normed spaces, we show that entropy and spherical entropy coincide for general convex bodies. 

In section \ref{sec_entropy_bounds}, we give the proofs of our main theorems. In the final section \ref{sec_centro_proj}, we give an intrinsic definition of the centro-projective surface area and study some of its properties. In particular, we show that it is upper semi-continuous with respect to Hausdorff topology.  
%--------------------------------------------------
\subsection*{Acknowledgements}
We wish to thank Bruno Colbois and Daniel Hug for interesting discussions and Franz Schuster for useful remarks on an earlier version of this paper. 
  
%------------------------------------------------
\section{Preliminaries on Convex bodies \\ and Hilbert Geometries}
\label{sec_preliminaries}

\subsection{Convex bodies} Let $V$ be a finite-dimensional real vector space. By \emph{convex body}, we mean a compact convex set $K\subset V$ with non-empty interior (note that this last condition is sometimes not required in the literature). Most of the time, the convex bodies will be assumed to contain the origin in their interiors. In such a case, we will call as usual \emph{Minkowski functional} the positive, homogeneous of degree one function whose level set at height 1 is the boundary $\partial K$. It is a convex function and by Alexandroff's theorem, it admits a quadratic approximation  almost everywhere (see e.g. A.D.~Alexandroff~\cite{alexandrov} or L.C.~Evans \& R.F. Gariepy~\cite[p. 242]{MTFPF92}). In the following, boundary-points where Alexandroff's theorem applies will be called \emph{smooth}. Assuming the vector space to be equipped with an inner product, the principal curvatures of the boundary and its Gauss curvature $k$ are well defined at every smooth point.

We will be concerned with generalizations and variations of \emph{Blaschke's rolling theorem}, a proof of which may be found in K.~Leichtwei\ss~\cite{leichtweiss93}.

\begin{theorem}[\bf W. Blaschke, \cite{blaschke56}]
Let $K$ be a convex body in $\dR^n$ whose boundary is $C^2$ with everywhere positive Gaussian curvature. Then there exist two positive radii $R_1$ and $R_2$ such that for every boundary point $p$, there exists a ball of radius $R_1$ (resp. $R_2$) containing $p$ on its boundary and contained in $K$ (resp. containing $K$).
\end{theorem}
We first remark that for the ``inner part'' of Blaschke's result, the regularity of the boundary may be lowered. Recall that the boundary of a convex body is $C^{1,1}$ provided it is $C^1$ and the Gauss map is Lipschitz-continuous. 
Roughly speaking, the second condition says that the curvature of the 
boundary remains bounded, even if it is only almost everywhere defined. The following proposition then gives a geometrical 
characterization of such bodies, see L.~H\"ormander~\cite[proposition~2.4.3]{hor94}
 or V.~Bangert~\cite{bang99} and D.~Hug~\cite{hug99}.

\begin{proposition}\label{c11}
The boundary of a convex body $K$ is $C^{1,1}$ if and only if there exists some $R>0$ such that $K$ is the union of balls with radius $R$.
\end{proposition}

Without assumption on the boundary, there is still an integral version of Blaschke's rolling theorem. 

\begin{theorem}[\bf C. Sch\"utt \& E. Werner, \cite{schwer90}] \label{thm_schuett_werner} 
 For a  convex body $K$ containing the unit ball of a Euclidean space and $p \in \partial K$, let $R(p) \in [0,\infty)$ be the radius of the biggest ball contained in $K$ and containing $p$. Then for all $0<\alpha<1$
\begin{equation}
 \int_{\partial K} R^{-\alpha} d\mathcal{H}^{n-1} < \infty.
\end{equation}
\end{theorem}

We will need the following refinement of this theorem.
\begin{proposition} \label{proposition_blaschke_strong}
In the same situation as in Theorem \ref{thm_schuett_werner}, 
for each Borel subset $B \subset \partial K$ we have 
\begin{multline}
 \int_B R^{-\alpha}d\mathcal{H}^{n-1} \leq\\ 2(n-1)^\alpha\left(\frac{2^\alpha}{1-2^{\alpha-1}}\right)^\alpha\left(\mathcal{H}^{n-1}(B)\right)^{1-\alpha} \left(\cH^{n-1}(\partial K)\right)^\alpha.
\end{multline}
In particular for some constant $C$ depending on $K$ we have
\begin{equation}
 \int_B R^{-\frac12}d\mathcal{H}^{n-1} \leq C \left(\mathcal{H}^{n-1}(B)\right)^{\frac12}\text{.}
\end{equation}
\end{proposition}

\proof
By (\cite{schwer90}, Lemma 4), we have for $0 \leq t \leq 1$
\begin{equation}\label{keyingredient}
 \mathcal{H}^{n-1}\bigl(\{p \in \partial K| R(p) \leq t\}\bigr)\leq (n-1)t\ \mathcal{H}^{n-1}(\partial K)\text{,} 
\end{equation}
from which we deduce that, for each $0<\epsilon<1$
\begin{multline}
 \int_{\partial K \cap \{R<\epsilon\}} R^{-\alpha} d\cH^{n-1} = \sum_{i=0}^\infty \int_{\partial K \cap \{\epsilon 2^{-i-1} \leq R<2^{-i} \epsilon\}} R^{-\alpha} d\cH^{n-1}\\
 \leq \sum_{i=0}^\infty (\epsilon 2^{-i-1})^{-\alpha}\ \cH^{n-1} \bigl(\partial K \cap \{\epsilon 2^{-i-1} \leq R<2^{-i} \epsilon\}\bigr)\\
 \leq \sum_{i=0}^\infty (\epsilon 2^{-i-1})^{-\alpha} (n-1)2^{-i}\epsilon\ \cH^{n-1} (\partial K)\\
 = \epsilon^{1-\alpha} (n-1) \frac{2^\alpha}{1-2^{\alpha-1}}\ \cH^{n-1} (\partial K).
\end{multline}

It follows that
\begin{align*}
 \int_B R^{-\alpha}d\mathcal{H}^{n-1} & = \int_{B \cap \{R<\epsilon\}} R^{-\alpha}d\mathcal{H}^{n-1}+\int_{B \cap \{R \geq \epsilon\}} R^{-\alpha}d\mathcal{H}^{n-1}\\
& \leq \epsilon^{1-\alpha} (n-1) \frac{2^\alpha}{1-2^{\alpha-1}}\ \cH^{n-1} (\partial K)+ \epsilon^{-\alpha}\ \mathcal{H}^{n-1}(B).
\end{align*}

Choosing 
\begin{displaymath}
 \epsilon:=\frac{1-2^{\alpha-1}}{2^\alpha(n-1)}\ \frac{\cH^{n-1} (B)}{\cH^{n-1}(\partial K)}
\end{displaymath}

yields the inequality of the lemma.
\endproof 

\subsection{Hilbert geometries} Given two distinct points $x, y \in \inte K$, the \emph{Hilbert distance} between $x$ and $y$ is defined by 
\begin{displaymath}
d(x,y):=\frac12 \bigl|\log [a,b,x,y]\bigr|,
\end{displaymath}
where $a$ and $b$ are the intersections of the line passing through $x$ and $y$ with the boundary $\partial K$, and $[a,b,x,y]$ denotes the cross-ratio (with the convention of \cite{bridson_haefliger}).

This distance is invariant under projective transformations. If $K$ is an ellipsoid, the Hilbert geometry on $\inte K$ is isometric to hyperbolic $n$-space. 

Unbounded closed convex sets with non-empty interiors and not containing a straight line are projectively equivalent to convex bodies. Therefore, the definition of the distance naturally extends to the interiors of such convex sets. In particular the convex sets bounded by parabolas are also isometric to the hyperbolic space.

Let us assume the origin $o$ lies inside the interior of $K$. We will write $B(r)$ for the \emph{metric ball} of radius $r$ and centered at $o$. Its boundary, the \emph{metric sphere}, will be denoted by $S(r)$. Let $a\colon\partial K\longrightarrow\dR_+$ be defined by the equation
\begin{displaymath}
-a(p)p\in\partial K, 
\end{displaymath}
so the letter $a$ refers to the antipodal point. It is an easy exercise to check that metric spheres are parameterized by the boundary $\partial K$ as
\begin{displaymath}
S(r)=\bigl\{\phi(p,r)\colon p\in\partial K\bigr\}, 
\end{displaymath}
where
\begin{align}\label{param_spheres}
\phi\colon \partial K\times\dR_+ & \to \inte{K}\\
(p,r) & \mapsto a\frac{e^{2r}-1}{a e^{2r}+1}\ p.\nonumber
\end{align}

The Hilbert distance comes from a Finsler metric on the interior of $K$. Given $x \in \inte K$ and $v \in T_xV$, the Finsler norm of $v$ is given by 
\begin{equation} \label{eq_finsler_norm}
	\|v\|_x=\frac12 \left(\frac{1}{t_1}+\frac{1}{t_2}\right),
\end{equation} 
where $t_1,t_2 >0 $ are such that $x\pm t_iv\in \partial K$. Again, we do not exclude that one of the $t_i$'s is infinite. Equivalently, if $F_x$ is the Minkowski functional of $K-x$, then 
\begin{displaymath}
\Vert v\Vert_x=\frac{1}{2}\bigl(F_x(v)+F_x(-v)\bigr). 
\end{displaymath}

The Finsler metric makes it possible to measure the length of a differentiable curve $c:I \to \inte K$ by
\begin{displaymath}
	l(c):=\int_I \bigl\|c'(t)\bigr\|_{c(t)} dt. 
\end{displaymath} 

It is less trivial to measure the area (or volume) of higher dimensional subsets of $\inte K$. In fact, different notions of volume are being used. The most important ones are the Busemann definition (which equals the Hausdorff $n$-dimensional measure) and the Holmes-Thompson definition. In the following, only the axioms of a \emph{volume} as defined in \cite{alth04} will be used. We will make use of the following properties:
\begin{itemize}
 \item $\vol$ is a Borel measure on $\inte K$ which is absolutely continuous with respect to Lebesgue measure.
\item If $A \subset K \subset L$, where $K,L$ are compact convex sets, then the measure of $A$ with respect to $K$ is larger than the measure of $A$ with respect to $L$.
\item If $K$ is an ellipsoid, then $\vol(A)$ is the hyperbolic volume of $A$.  
\end{itemize}

The following projective invariants of convex bodies will be our main subjects of investigation.
\begin{definition}
The \emph{upper (resp. lower) volume entropy} of $K$ is 
\begin{align*}
\overline{\ent}(K) & :=\limsup_{r\to\infty}\frac{\log\bigl(\vol B(r)\bigr)}{r}, \\
\underline{\ent}(K) & :=\liminf_{r\to\infty}\frac{\log\bigl(\vol B(r)\bigr)}{r}. 
\end{align*}
If the upper and lower volume entropies of $K$ coincide, their common value is called volume entropy of $K$ and denoted by $\ent K$. 
\end{definition}
Note that these invariants are independent of the choice of the center and of the choice of the volume definition.

\subsection{Busemann's density}

For simplicity, we restrict ourselves to Busemann's volume, although all results remain true for every other choice of volume. The reason is that the proofs of the crucial propositions \ref{proposition_parabola} and \ref{prop_pointwise_limit} below do not use any particular property of Busemann's volume, but only the axioms satisfied by every definition of volume.  

The density of Busemann's volume (with respect to some Lebesgue measure $\cL$) is given by 
\begin{displaymath}
 \sigma(x)=\frac{\omega_n}{\cL(B_x)},
\end{displaymath}
where $B_x$ is the tangent unit ball of the Finsler metric at $x$ and $\omega_n$ is the (Euclidean) volume of the unit ball in $\dR^n$. The volume of a Borel subset $A \subset \inte K$ is thus given by
\begin{displaymath}
 \vol(A)=\int_A \sigma \ d\cL. 
\end{displaymath}

We now state and prove some propositions concerning upper bounds and asymptotic behaviors of Busemann's densities for points which are close to the boundary of particular convex sets. We will make use of an auxiliary inner product, calling $\cL$ and $\mu$ the corresponding Lebesgue measure and volume $n$-form. Busemann densities are defined with this particular choice of measure.

\begin{proposition} \label{proposition_truncation}
 Let $K,K'$ be closed convex sets not containing any straight line and $\sigma:\inte K \to \dR$,  $\sigma':\inte K' \to \dR$ their corresponding Busemann densities.
Let $p \in \partial K$, $E_0$ a support hyperplane of $K$ at $p$ and $E_1$ a hyperplane parallel to $E_0$ intersecting $K$. Suppose that $K$ and $K'$ have the same intersection with the strip between $E_0$ and $E_1$ (in particular $p \in \partial K'$). Then 
\begin{displaymath}
 \lim_{y \to p} \frac{\sigma(y)}{\sigma'(y)}=1.
\end{displaymath}
\end{proposition}
\proof
Let $d$ be the distance between $E_0$ and $E_1$ and $(y_i)$ a sequence of points of $\inte K$ converging to $p$. We may suppose that the distance $d_i$ between $y_i$ and $E_0$ is strictly less than $d$. For every fixed point $y_i$ and non-zero tangent vector $v \in T_{y_i}K$, let $t_1,t_2 \in \dR_+ \cup \{\infty\}$ be such that $y_i \pm t_{1,2}v \in \partial K$; let $t_1',t_2'$ be the corresponding numbers for $K'$. Since at least one of $y_i+t_1v$ and $y_i-t_2v$ is inside the strip, say $y_i+t_1v$, we must have $t_1=t_1'$. 

Either $t_2=t_2'$ and $\|v\|_i=\|v\|_i'$, or $t_2 \neq t_2'$, in which case 
\begin{displaymath}
\frac{t_1}{t_2}, \frac{t_1'}{t_2'}\leq \frac{d_i}{d-d_i}. 
\end{displaymath}

Therefore, 
\begin{displaymath}
 \frac{d-d_i}{d}\leq\frac{\|v\|_i}{\|v\|_i'}=\frac{1+\frac{t_1}{t_2}}{1+\frac{t_1'}{t_2'}}\leq\frac{d}{d-d_i}
\end{displaymath}
which shows that, as functions on $\dR P^{n-1}$, $\Vert\cdot\Vert_i/\Vert\cdot\Vert_i'$ uniformly converge to 1.
Hence, for every $\epsilon$ and every $i$ large enough, 
\begin{displaymath}
(1-\epsilon)B_{y_i}\subset B_{y_i}'\subset(1+\epsilon)B_{y_i}, 
\end{displaymath}
which implies the convergence of $\sigma/\sigma'$ to 1.
\endproof

\begin{proposition} \label{proposition_parabola}
 Let $V=\dR^n$ with its usual scalar product. Let $P$ be the convex set bounded by the parabola $y=\sum_{i=1}^{n-1} \frac{c_i}{2} x_i^2, c_1,\ldots,c_{n-1}>0$. Then 
\begin{equation}
 \sigma(0,\ldots,0,1-\lambda)=\frac{\sqrt{c}}{\bigl(2(1-\lambda)\bigr)^{\frac{n+1}{2}}},
\end{equation}
where $c=\prod_{i=1}^{n-1} c_i$. 
\end{proposition}

\proof
By the invariance of the Hilbert metric under projective transformations, the tangent unit sphere at any point of $\inte P$ is an ellipse. At the point $(0,\ldots,0,1-\lambda)$, the symmetry implies that the principal axes of this ellipse are parallel to the coordinate axes. Hence 
\begin{displaymath}
 \sigma=\frac{1}{\prod_{i=1}^n l_i},
\end{displaymath}
where the $l_i$'s, $i=1,\ldots,n$, are the Euclidean lengths of the principal half-axes. 

Now $l_i=\sqrt{\frac{2(1-\lambda)}{c_i}}, i=1,\ldots,n-1$ and $l_n=2(1-\lambda)$. 
\endproof

\begin{proposition} \label{prop_pointwise_limit}
Assume the origin $o$ is inside $\inte{K}$. For a smooth point $p$ of $\partial K$, let $n(p)$ be the outward normal vector and let $k(p)$ be the Gauss curvature of $\partial K$ at $p$. Then   
\begin{equation}
 \lim_{\lambda \to 1} \sigma(\lambda p) (1-\lambda)^\frac{n+1}{2}=\frac{\sqrt{k(p)}}{\Bigl(2\bigl\langle p,n(p)\bigr\rangle\Bigr)^\frac{n+1}{2}}.
\end{equation}
\end{proposition}

\proof
Let us choose a frame $(p;v_1,\ldots,v_{n-1},v_n)$ where $v_1,\ldots,v_{n-1} \in T_p \partial K$ are unit vectors tangent to the principal curvature directions of $\partial K$ at $p$ and $v_n=-p$.  In these coordinates, the boundary of $K$ is locally the graph of a function:  $y=\sum_{i=1}^{n-1} \frac{c_i}{2} x_i^2+R(|x|)$ with $R(|x|)=o(|x|^2)$ and $c_1,\ldots,c_{n-1} \geq 0$. We set
\begin{displaymath}
 c:=\prod_{i=1}^{n-1} c_i.
\end{displaymath}
A small computation shows that 
\begin{displaymath}
dx_1 \wedge \ldots \wedge dx_{n-1} \wedge dy=\frac{1}{m} \mu,
\end{displaymath}
where $\mu$ is the Euclidean $n$-form and $m:=\mu(v_1,\ldots,v_n)=\bigl\langle p,n(p)\bigr\rangle$. Also, the Gauss curvature at $p$ is given by 
\begin{displaymath}
k(p)=cm^{n-1}.
\end{displaymath}

Let us fix $\epsilon>0$. Locally, the parabola defined by 
\begin{displaymath}
y=\sum_{i=1}^{n-1} \frac{c_i+\epsilon}{2}x_i^2
\end{displaymath}
lies inside $K$. Cutting it with some horizontal hyperplane, we obtain a convex body $K'$ inside $K$. In particular, the metric of $K'$ is greater than or equal to the metric of $K$, hence $\sigma'(\lambda p) \geq \sigma(\lambda p)$ for $\lambda$ near $1$.

Then by propositions \ref{proposition_truncation} and \ref{proposition_parabola}, 
\begin{align} \label{eq_limsup}
\nonumber \limsup_{\lambda \to 1} \sigma(\lambda p)(1-\lambda)^\frac{n+1}{2} & \leq \lim_{\lambda \to 1} \sigma'(\lambda p)(1-\lambda)^\frac{n+1}{2}\\
& = \frac{\sqrt{\prod_{i=1}^{n-1} (c_i+\epsilon)}}{2^\frac{n+1}{2}m}.   
\end{align}
Note that $\sigma>0$, hence this already settles the case $k=c=0$ since $\epsilon$ was arbitrary small.

If $c>0$ and $0<\epsilon<\min\{c_1,\ldots,c_{n-1}\}$, the parabola $P$ defined by 
\begin{displaymath}
 y=\sum_{i=1}^{n-1} \frac{c_i-\epsilon}{2}x_i^2
\end{displaymath}
locally contains $K$. Cutting it with some horizontal hyperplane, we obtain a convex body $K'$ inside $P$. By propositions \ref{proposition_truncation} and \ref{proposition_parabola} again, 
\begin{align} 
\nonumber \liminf_{\lambda \to 1} \sigma(\lambda p)(1-\lambda)^\frac{n+1}{2} &  \geq \liminf_{\lambda \to 1} \sigma'(\lambda p)(1-\lambda)^\frac{n+1}{2} \\
& = \frac{\sqrt{\prod_{i=1}^{n-1} (c_i-\epsilon)}}{2^\frac{n+1}{2}m}. \label{eq_liminf}
\end{align}

From \eqref{eq_limsup} and \eqref{eq_liminf} (with $\epsilon \to 0$) we get 
\begin{displaymath}
 \lim_{\lambda \to 1} \sigma(\lambda p)(1-\lambda)^\frac{n+1}{2} = \frac{\sqrt{c}}{2^\frac{n+1}{2}m}.
\end{displaymath}
\endproof

To state precisely our main theorem in section \ref{sec_entropy_bounds} we need to introduce the pseudo-Gauss curvature of the boundary of a convex set $K$ in $\dR^n$.

For a smooth point $p\in \partial K$, let $n(p)$ be the outward normal of $\partial K$ at $p$. For each unit vector $e \in T_p \partial K$, let $H_e(p)$ be the affine plane containing $p$ and directed by the
vectors $e$ and $n(p)$. We define $R_e$ as the radius of 
the biggest disc containing $p$ inside $K_e:=K \cap H_e(p)$.

\begin{definition}\label{pseudogauss}
The \textit{pseudo Gauss-curvature} $\bar k(p)$ of $\partial K$ at $p$ is the minimum of the numbers 
\begin{displaymath}
 \prod_{i=1}^{n-1} R_{e_i}(p)^{-1},
\end{displaymath}
where $e_1,\ldots,e_{n-1}$ ranges over all orthonormal bases of $T_p \partial K$.  
\end{definition}

\begin{proposition}\label{proposition_blaschke_trick}
  Let $V$ be a Euclidean vector space of dimension $n$. Let $K$ be a convex body
containing the unit ball $B$. Then for $\frac12 \leq \lambda < 1$ and $p \in \partial K$ 
\begin{equation}
\sigma(\lambda p) \leq \frac{\omega_n n!}{2^n(1-\lambda)^\frac{n+1}{2}} \bar k(p)^{1/2}\text{.}
\end{equation}
\end{proposition}

\proof
We use the same notation as in the definition of $\bar k$.   
We may suppose that for all $i$, $R_i:=R_{e_i}(p)>0$, 
otherwise the statement is trivial. By definition of $R_i$, there is a $2$-disc
 $B_i(p)$ of radius $R_i$ inside $K_{e_i}$ containing $p$. Let us denote by $B(e_i)$ the intersection of
$B$ with the affine plane $p+H_{e_i}$.
Since $B(e_i),B_i(p)\subset K$, one has
\begin{displaymath}
\hat C_i:=conv\left(B(e_i)\times\{0\}\cup B_i(p)\times\{1\}\right)\subset K_{e_i}\times[0,1]. 
\end{displaymath}

Note that $\hat C_i$ is a truncated cone. Let $E_i$ be the plane containing the line that is parallel to $T_p\partial K_{e_i}$ and that passes through the points $o\times\{0\}$ and $p\times\{1\} $. 
With $\pi:V \times [0,1] \to V$ the projection on the first component, $C_i:=\pi(E_i\cap\hat C_i)\subset K$ is bounded by a truncated conic.

In the non-orthogonal frame $(o;p,e_i)$, $C_i$ is given by
\begin{displaymath}
(2R_i-1)x^2+2(1-R_i)x+ y_1^2 \leq 1,\quad 0 \leq x\leq 1. 
\end{displaymath}

Now let $C$ be the convex hull of the union of the $C_i$.
Then the  polytope $P$ with vertices
\begin{displaymath}
 \left(\lambda, 0,\ldots,\pm \sqrt{(1-\lambda)(2\lambda R_i-\lambda+1)},0,\ldots,0\right), (1,\vec0), (2\lambda-1,\vec0)
\end{displaymath}
lies inside $C$, with all but the last vertex being on the boundaries of the respective $C_i$'s. 

Its volume is given by
\begin{align}
\nonumber \cL(P) & =\frac{2^n \bigl\langle p,n(p)\bigr\rangle }{n!}(1-\lambda)^\frac{n+1}{2}\prod_{i=1}^{n-1}(2\lambda R_i-\lambda+1)^\frac{1}{2}\\
\nonumber & \geq \frac{2^n}{n!}(1-\lambda)^\frac{n+1}{2} (R_1\cdot R_2\cdots R_{n-1})^{\frac{1}{2}}\\
& = \frac{2^n}{n!}(1-\lambda)^\frac{n+1}{2} \bar k^{-\frac{1}{2}}(p).
\end{align}
The factor $\bigl\langle p,n(p)\bigr\rangle$ in the first line is due to the fact that our coordinate system is not orthonormal. Since the unit ball is contained in $K$, this factor is at least $1$. 

From $P \subset C \subset K$ and the fact that $P$ is centered at $\lambda p$, we deduce that 
\begin{displaymath}
 \sigma(\lambda p) \leq \frac{\omega_n}{\cL(P)} \leq \frac{\omega_n n!}{2^n} (1-\lambda)^{-\frac{n+1}{2}} \bar k^{\frac{1}{2}}(p).
\end{displaymath}
\endproof

The next proposition will be needed in the construction of a convex body with entropy between $0$ and $1$. 
 
\begin{proposition} \label{proposition_density_triangle}
 Let $K=oab$ be a triangle with $1 \leq oa,ob \leq 2$ and such that the distance from $o$ to the line passing through $a$ and $b$ is at least $1$. 
 Let $p$ be a point in the interior of the side $ab$ and suppose that $\min\{ap,bp\}\geq \epsilon >0$. Then for $\lambda \geq \frac12$ Busemann's density of $K$ at $\lambda p$ is bounded above by 
\begin{displaymath}
\sigma(\lambda p) \leq 32 \pi\max\left\{\frac{1}{\epsilon (1-\lambda)},\frac{1}{\epsilon^2}\right\}.
\end{displaymath}
\end{proposition}

\proof
The hypothesis on the triangle implies that $\sin(abo),\sin(bao) \geq \frac12$.  

Let $a'$ be the intersection of the line passing through $a$ and $z:=\lambda p$ with $ob$ and define $b'$ similarly. 

The unit tangent ball at $z$ is a hexagon centered at $z$. The length of one of its  half-diagonals is the harmonic mean of $za$ and $za'$; the length of the second half-diagonal is the harmonic mean of $zb$ and $zb'$ and the third half-diagonal has length $\frac{2 op}{\frac{1}{\lambda}+\frac{1}{1-\lambda}} \geq 1-\lambda$.  

An easy geometric argument shows that $za',zb \geq \frac12 pb \sin(abo) \geq \frac14 \epsilon$ and $za,zb' \geq \frac12 pa \sin(bao) \geq \frac14 \epsilon$. 

The area $A$ of the hexagon is at least half of the minimal product of two of its half-diagonals, hence 

\begin{displaymath}
 A \geq \min\left\{\frac18 \epsilon (1-\lambda),\frac{1}{32}\epsilon^2\right\}.
\end{displaymath}
\endproof

\subsection{Volume entropy of spheres}
By definition, the entropy controls the volume growth of metric balls in Hilbert geometries. We show in this section that it coincides with the growth of areas of metric spheres. Again, there are several definitions of area of hypersurfaces in Finsler geometry. For simplicity, we consider Busemann's definition which gives the Hausdorff $(n-1)$-measure of these hypersurfaces.

We will need the following two lemmas:
\begin{lemma}[Rough monotonicity of area]
There exist a monotone function $f$ and a constant $C_1>1$ such that for all $r>0$
\begin{equation} \label{eq_comparison_holmes_thompson}
 C_1^{-1} f(r) \leq Area(S(r)) \leq C_1 f(r).
\end{equation}
\end{lemma}

\proof
Let $f(r)$ be the Holmes-Thompson area of $S(r)$. Since all area definitions agree up to some universal constant, inequality \eqref{eq_comparison_holmes_thompson} is trivial. It remains to show that $f$ is monotone.

If $\partial K$ is $C^2$ with everywhere positive Gaussian curvature then the tangent unit spheres of the Finsler metric are quadratically convex. According to \cite[theorem~1.1 and remark~2]{alfe98} there exists a Crofton formula for the Holmes-Thompson area, from which the monotonicity of $f$ easily follows. 

Such smooth convex bodies are dense in the set of all convex bodies for the Hausdorff topology (see e.g. \cite[lemma~2.3.2]{hor94}). By approximation, it follows that $f$ is monotone for arbitrary $K$.
\endproof

\begin{lemma}[Co-area inequalities]
There exists a constant $C_2>1$ such that for all $r>0$
\begin{displaymath}
 C_2^{-1} Area(S(r)) \leq \frac{\partial}{\partial r}\vol(B(r)) \leq C_2 Area(S(r)).
\end{displaymath}
\end{lemma}

\proof
Let $\mu:=\sigma dx_1\land\dots\land dx_n$ be the volume form, and let $\alpha$ be the $n-1$-form on $S(r)$ whose integral equals the area.  

Since
\begin{displaymath}
\vol(B(r))=\int_0^r\int_{S(s)}i_{\partial_r}\mu\ ds, 
\end{displaymath}
where $\partial_r$ at $\lambda p\in S(s)$ is the tangent vector multiple of $\vec{op}$ with unit Finsler norm, we have to compare $i_{\partial_r}\mu$ and $\alpha$. 

We will assume that $S(r)$ is differentiable at $\lambda p$. The section of the unit tangent ball by the tangent space $T_{\lambda p}S(r)$ will be called $\gamma$. By definition of Busemann area, the area of $\gamma$ measured with the form $\alpha$ is the constant
\begin{displaymath}
\alpha(\gamma)=\omega_{n-1}. 
\end{displaymath}

In the same way, calling $\Gamma$ the half unit ball containing $\partial_r$ and bounded by $\gamma$, one has 
\begin{displaymath}
\mu(\Gamma)=\frac{1}{2}\omega_n. 
\end{displaymath}

Since $\Gamma$ is convex it contains the cone with base $\gamma$ and vertex $\partial_r$. Therefore,
\begin{equation}\label{upper}
\frac{1}{n}i_{\partial_r}\mu(\gamma)\leq\frac{1}{2}\omega_n.
\end{equation}

By Brunn's theorem (see e.g. \cite[theorem 2.3]{kold2005}), the sections of the tangent unit ball with hyperplanes parallel to $\gamma$ have an area lesser than or equal to the area of $\gamma$. Also the tangent unit ball has a supporting hyperplane at $\partial_r$ which is parallel to $\gamma$. Therefore, by Fubini's theorem, the cylinder $\gamma\times([0,1]\cdot\partial_r)$ has a volume greater than or equal to the volume of $\Gamma$ (even if it generally does not contain $\Gamma$). Hence,
\begin{equation}\label{lower}
\frac{1}{2}\omega_n\leq i_{\partial_r}\mu(\gamma).
\end{equation}

Inequalities ~\eqref{upper} and \eqref{lower} give
\begin{displaymath}
\frac{1}{2}\frac{\omega_n}{\omega_{n-1}}\alpha(\gamma)\leq i_{\partial_r}\mu(\gamma)\leq\frac{n}{2}\frac{\omega_n}{\omega_{n-1}}\alpha(\gamma), 
\end{displaymath}
from which the result easily follows.
\endproof

\begin{theorem}
 The spherical entropy coincides with the entropy. More precisely,
\begin{align*}
\limsup_{r \to \infty} \frac{\log \ar(S(r))}{r} & =\overline{\ent}K,\\
 \liminf_{r \to \infty} \frac{\log \ar(S(r))}{r} & =\underline{\ent}K.
\end{align*}
\end{theorem}

\proof
For convenience, let 
\begin{align*}
 V(r) & :=\vol B(r),\\ 
 A(r) & :=\ar S(r).
\end{align*}

Using the previous two lemmas, one has for all $r>0$
\begin{multline*}
 V(r)=\int_0^r V'(s)ds \leq C_2 \int_0^r A(s)ds \leq C_1C_2\int_0^r f(s)ds\\ \leq C_1C_2 f(r)r \leq C_1^2C_2 A(r)r.
\end{multline*}

It follows that 
\begin{multline*}
 \overline{\ent}K=\limsup_{r \to \infty} \frac{\log V(r)}{r} \leq \limsup_{r \to \infty} \frac{\log C_1^2C_2A(r)r}{r}\\=\limsup_{r \to \infty} \frac{\log \ar(S(r))}{r}. 
\end{multline*}

Similarly, for each $\epsilon>0$ 
\begin{multline*}
 V(r(1+\epsilon)) = \int_0^{r(1+\epsilon)} V'(s)ds \geq C_1^{-1}C_2^{-1}\int_0^{r(1+\epsilon)} f(s)ds\\ \geq C_1^{-1}C_2^{-1} \int_r^{r(1+\epsilon)} f(s)ds  \geq C_1^{-1}C_2^{-1} f(r) r \epsilon \geq C_1^{-2}C_2^{-1} A(r) r \epsilon  
\end{multline*}
and hence
\begin{multline*}
 (1+\epsilon)\overline{\ent} K = (1+\epsilon) \limsup_{r \to \infty} \frac{\log V(r(1+\epsilon))}{r(1+\epsilon)} \geq \limsup_{r \to \infty} \frac{\log C_2^{-1}C_1^{-2} A(r)r\epsilon}{r} \\
= \limsup_{r \to \infty} \frac{\log \ar(S(r))}{r}. 
\end{multline*}

Letting $\epsilon \to 0$ gives the first equality. The second one follows in a similar way. 
\endproof

%---------------------------------------------------
\section{Entropy bounds}
\label{sec_entropy_bounds}

\subsection{Upper entropy bound in arbitrary dimension}

We may now state and prove our first main theorem. 

\begin{theorem} \label{thm_main_thm}
Let $K$ be an $n$-dimensional convex body and $o \in \inte K$. For any point $p \in \partial K$ we denote by  
$\bar k(p)$ its pseudo-Gauss curvature as in definition  \ref{pseudogauss}. If 
\begin{equation} \label{eq_blaschke_hypothesis}
 \int_{\partial K} {\bar k}^{\frac12}(p)dp<\infty,
\end{equation}
then 
\begin{equation} \label{eq_entropy_coefficient}
 \lim_{r \to \infty} \frac{\vol B(o,r)}{\sinh^{n-1}r}=\frac{1}{n-1} \cA_p(K). 
\end{equation}
In particular, 
\begin{displaymath}
 \overline{\ent} K \leq n-1,
\end{displaymath}
and if $\cA_p(K) \neq 0$, then $\overline{\ent} K=n-1$.
\end{theorem}

\proof
Using the parameterization~\eqref{param_spheres}, the volume of metric balls is given by
\begin{displaymath}
\vol(B(r))=\int_0^r\int_{\partial K}F(p,r)\ d\cH^{n-1}, 
\end{displaymath}
where
\begin{displaymath}
F(p,r):=\sigma \bigl(\phi(p,r)\bigr)\Jac \phi(p,r). 
\end{displaymath}

The Jacobian may be explicitly computed: 
\begin{displaymath}
 \Jac \phi(p,r)=\frac{(e^{2r}-1)^{n-1} e^{2r}}{(ae^{2r}+1)^{n+1}} 2a^n (1+a) \bigl\langle p,n(p)\bigr\rangle. 
\end{displaymath}
In particular,
\begin{equation} \label{eq_limit_jac}
 \lim_{r \to \infty} e^{2r} \Jac \phi(p,r)=\frac{2(1+a)\bigl\langle p,n(p)\bigr\rangle}{a}.
\end{equation}

On the other hand, for each smooth boundary point $p$ we have, by proposition \ref{prop_pointwise_limit}, 
\begin{equation} \label{eq_limit_kappa}
 \lim_{r \to \infty} \frac{\sigma\bigl(\phi(p,r)\bigr)}{e^{(n+1)r}}=\frac{\sqrt{k(p)}}{\Bigl(2\bigl\langle p,n(p)\bigr\rangle\Bigr)^\frac{n+1}{2}} \frac{a^\frac{n+1}{2}}{(1+a)^\frac{n+1}{2}}.
\end{equation}

Then, by proposition \ref{proposition_blaschke_trick} 
and the hypothesis \eqref{eq_blaschke_hypothesis},
\begin{align}\label{eq_change_order}
\lim_{r\to\infty}\frac{1}{e^{(n-1)r}}\int_{\partial K}F(p,r)d\cH^{n-1}
& =\int_{\partial K}\lim_{r\to\infty}\frac{F(p,r)}{e^{(n-1)r}}d\cH^{n-1}\\
& =\int_{\partial K}\lim_{r\to\infty}\frac{\sigma \bigl(\phi(p,r)\bigr)}{e^{(n+1)r}} \lim_{r \to \infty} e^{2r} \Jac \phi(p,r) d\cH^{n-1}\\
&= \int_{\partial K} \frac{\sqrt{k(p)}}{\Bigl(2\bigl\langle p,n(p)\bigr\rangle\Bigr)^\frac{n-1}{2}} \left(\frac{a}{1+a}\right)^\frac{n-1}{2} d\cH^{n-1}\nonumber\\
&= \frac{1}{2^{n-1}}\cA_p(K).\nonumber
\end{align}

By L'Hospital's rule we get
\begin{displaymath}
\lim_{r \to \infty} \frac{\vol\bigl(B(r)\bigr)}{e^{(n-1)r}}=\lim_{r\to\infty}\dfrac{\int_0^r\int_{\partial K}F(p,s)d\cH^{n-1}ds}{(n-1)\int_0^r e^{(n-1)s}ds}=\dfrac{1}{2^{n-1}(n-1)}\cA_p(K).
\end{displaymath}
\endproof

{\bf Remark:} The metric balls $B(r)$ are projective invariants of $K$. There is an affine version of the previous theorem using the affine balls $B_a(r):=\tanh(r)K$ (where multiplication is with respect to the center $o$). Under the same assumptions as in theorem \ref{thm_main_thm}, we obtain that 
\begin{displaymath}
  \lim_{r \to \infty} \frac{\vol B_a(r)}{e^{(n-1)r}}=\frac{1}{2^{n-1}(n-1)} \cA_a(K)
\end{displaymath}
where $\cA_a(K)$ is the centro-affine area (see section~\ref{sec_centro_proj}). 
The proof goes as the previous one by replacing the function $a$ by $1$.

\begin{corollary}\label{corollary1} 
 Suppose $K$ is an $n$-dimensional convex body of class $C^{1,1}$. Then 
\begin{displaymath}
 \ent K=n-1.
\end{displaymath}
\end{corollary}

\proof
For any $p\in\partial K$, $R(p)$ is the biggest radius of a ball in $K$ containing $p$. 
By proposition~\ref{c11}, there exists a constant $R>0$ such that $R(p) \geq R$ for all $p \in \partial K$. 
It follows that the hypothesis \eqref{eq_blaschke_hypothesis} is satisfied and therefore $\ent K \leq n-1$. 

The Gauss map $\cG\colon\partial K\rightarrow S^{n-1}$ is well-defined and continuous. As a consequence of theorem~2.3 in Hug~\cite{hugII} and equation~2.7 in Hug~\cite{hugI}, the standard measure on the unit sphere is the push-forward of $k\cdot d\cH^{n-1}$, i.e.
\begin{displaymath}
\cG_*(k\cdot d\cH^{n-1}\mbox{}_{|\partial K})=d\cH^{n-1}|_{S^{n-1}}\ ,
\end{displaymath}
hence the curvature has a positive integral. Therefore, $\cA_p(K)>0$, and equation~\eqref{eq_entropy_coefficient} implies that $\ent K=n-1$. 
\endproof

\begin{corollary}
 If $K$ is an arbitrary $n$-dimensional convex body with $\cA_p(K) \neq 0$, then $\underline{\ent}K \geq n-1$.
\end{corollary}

\proof
Arguing as in the proof of theorem \ref{thm_main_thm} and using Fatou's lemma instead of the dominated convergence theorem gives the result.  
\endproof

\subsection{The plane case} 

Let us now assume that $n=2$. By theorem \ref{thm_schuett_werner}, the hypothesis \eqref{eq_blaschke_hypothesis} is satisfied for each convex body $K$. Therefore 

\begin{equation} \label{eq_easy_bound_plane}
 \overline{\ent} K \leq 1
\end{equation}
and 
\begin{displaymath}
 \lim_{r \to \infty} \frac{\vol B(o,r)}{\sinh r}=\cA_p(K).
\end{displaymath}

Next, we are going to prove a better bound for $\overline{\ent} K$. In order to state our main result, we need to recall some basic notions of measure theory in a Euclidean space and refer to P.~Mattila~\cite{mat99} for details. For a non-empty bounded set $A$, let $N(A,\epsilon)$ be the minimal number of $\epsilon$-balls needed to cover $A$. Then the upper Minkowski dimension of $A$ is defined as 
\begin{displaymath}
 \overline{\dim} A:=\inf\left\{s:\limsup_{\epsilon \to 0} N(A,\epsilon)\epsilon^s=0\right\}.
\end{displaymath}

One should note that this dimension is invariant under bi-Lipschitz maps. In particular, it does not depend on a particular choice of inner product and moreover it is invariant under projective maps provided the considered subsets are bounded.

Recall that a point $p\in K$ is called \emph{extremal} if it is not a convex combination of other points of $K$. The set of extremal points is a subset of $\partial K$, which we denote by $\ex K$. 

\begin{theorem} \label{thm_minkowski_bound}
Let $K$ be a plane convex body and $d$ be the upper Minkowski dimension of $\ex K$. Then the entropy of $K$ is bounded by 
\begin{displaymath}
 \overline{\ent K} \leq \frac{2}{3-d} \leq 1.
\end{displaymath}
\end{theorem}

\proof
Since the entropy is independent of the choice of the center, we may suppose that the Euclidean unit ball around $o$ is the maximum volume ellipsoid inside $K$. Then $K$ is contained in the ball of radius $2$ (see \cite{barv02}).  

Set $\epsilon:=e^{-\alpha r}$, where $\alpha \leq 1$ will be fixed later. Divide the boundary of $K$ into two parts:
\begin{displaymath}
 \partial K= \cB \cup \cG,
\end{displaymath}
where $\cB$ (the bad part) is the closed $\epsilon$-neighborhood around the set of extremal points of $K$ and $\cG$ (the good part) is its complement.

Using proposition~\ref{proposition_blaschke_strong} and equalities \eqref{eq_limit_jac}, \eqref{eq_limit_kappa}, we get the following upper bound for large values of $r$,
\begin{equation}
\int_{\frac{r}{2}}^r \int_{\cB} \sigma \bigl(\phi(p,s)\bigr) \Jac \phi(p,s) d\cH^1 ds 
 \leq O\left(e^r \sqrt{\cH^1(\cB)}\right). \label{eq_bad_part}
\end{equation}
  
Next, let $p \in \cG$. The endpoints of the maximal segment in $\partial K$ containing  $p$ are extremal points of $K$ and hence of distance at least $\epsilon$ from $p$. Therefore $K$ contains a triangle as in proposition~\ref{proposition_density_triangle} and if $s \geq r/2$ and $r$ is sufficiently large
\begin{displaymath}
 \sigma(\phi(p,s))=\sigma(\lambda \cdot p) \leq 32 \max\left\{\frac{1}{\epsilon (1-\lambda)},\frac{1}{\epsilon^2}\right\} =\frac{32}{\epsilon (1-\lambda)}.
\end{displaymath}

Integrating this from $r/2$ to $r$ yields
\begin{equation}
\int_\frac{r}{2}^r \int_{\cG} \sigma \bigl(\phi(p,s)\bigr)\Jac \phi(p,s) d\cH^1 ds = O\left(e^{\alpha r}\right). \label{eq_good_part}
\end{equation}

Let $d$ be the upper Minkowski dimension of the set of extremal points of $K$. Then, for each $\eta>0$, $N(\ex K,\epsilon)=o(\epsilon^{-d-\eta})$ as $\epsilon \to 0$. By definition of $N$, there is a covering of $\ex K$ by $N(\ex K,\epsilon)$ balls of radius $\epsilon$. Hence there is a covering of $\cB$ by $N(\ex K,\epsilon)$ balls of radius $2\epsilon$. The intersection of a $2\epsilon$-ball with $\partial K$ has length less than $4\pi \epsilon$. It follows that 
\begin{displaymath}
 \mathcal{H}^1(\cB) = o(\epsilon^{-d-\eta+1}).
\end{displaymath}

Since the volume of $B(r/2)$ is bounded by $O(e^{r/2})$ (see \eqref{eq_easy_bound_plane}), the volume of $B(r)$ is bounded by 
\begin{align*}
\vol B(r) & =\vol B(r/2)+\int_{\frac{r}{2}}^r \int_{\cB} \sigma \bigl(\phi(p,s)\bigr) \Jac \phi(p,s) d\cH^1 ds \\
& \quad +\int_\frac{r}{2}^r \int_{\cG} \sigma \bigl(\phi(p,s)\bigr)\Jac \phi(p,s) d\cH^1 ds \\ 
& = O(e^\frac{r}{2}) + O\bigl(e^{r (1-\frac{\alpha (1-d-\eta)}{2})}\bigr) + O\left(e^{\alpha r}\right).
\end{align*}

We fix $\alpha$ such that $1-\alpha \frac{1-d-\eta}{2}=\alpha$, i.e. $\alpha:=\frac{2}{3-d-\eta} >\frac23$. Then 
\begin{displaymath}
 \vol B(r)=O(e^{\alpha r}),
\end{displaymath}
which implies that the (upper) entropy of $K$ is bounded by $\alpha$. Since $\eta>0$ was arbitrary, the result follows. 
\endproof

\subsection{An example of non-integer entropy}

We will construct an example of a plane convex body with piecewise affine boundary whose entropy is strictly between $0$ and $1$. 

Let us choose a real number $s>2$ and set $\alpha_i:=\frac{C_s}{i^s}$ where $C_s>0$ is sufficiently small such that 
\begin{equation*} \label{eq_sum_angles}
3\sum_{i=1}^\infty \alpha_i < \pi.
\end{equation*}
Consider a centrally symmetric sequence $E$ of points on $S^1$ such that the angles between consecutive points are $\alpha_1,\alpha_1,\alpha_1,\alpha_2,\alpha_2,\alpha_2,\ldots$ (each angle appearing three times). 

\begin{theorem}
The entropy of $K=conv(E)$ is bounded by 
\begin{displaymath}
0<\frac{1}{s} \leq \underline{\ent} K \leq \overline{\ent} K \leq \frac{2s-2}{3s-4}<1.
\end{displaymath}
\end{theorem}

\proof
{\bf Lower bound}\\
The unit sphere of radius $r$ in the Hilbert geometry $K$ is $\tanh r K$ and consists of an infinite number of segments. 

An easy geometric computation shows that the middle segment $S_i(r)$ corresponding to $\alpha:=\alpha_i$ has for each $r \geq 0$ length bounded from below by 
\begin{displaymath}
l\bigl(S_i(r)\bigr) \geq \log\left(\frac{\tanh r}{1-\tanh r} \frac{2\sin \alpha/2 \sin(2\alpha)}{\cos \alpha/2} +1\right).
\end{displaymath}

Set 
\begin{displaymath}
 i_0(r):=\left\lfloor (2C_s)^{\frac{1}{s}} e^\frac{r}{s} \right\rfloor.
\end{displaymath}
Then, for sufficiently large $r$,  
\begin{displaymath}
\frac{\tanh r}{1-\tanh r} \frac{2\sin \alpha_i/2 \sin(2\alpha_i)}{\cos \alpha_i/2} \leq 1 \quad \forall i \geq i_0(r).
\end{displaymath}

By concavity of the $\log$-function, we have $\log(1+x) \geq x \log 2 \geq \frac{x}{2}$ for $0 \leq x \leq 1$. Therefore  
\begin{displaymath}
l\bigl(S(r)\bigr) \geq \frac12 \sum_{i=i_0}^\infty \frac{\tanh r}{1-\tanh r} \frac{2\sin \alpha_i/2 \sin(2\alpha_i)}{\cos \alpha_i/2}.
\end{displaymath}
For sufficiently large $r$, the first factor is bounded from below by $\frac{e^{2r}}{4}$, while the second is bounded from below by $\alpha_i^2$. We thus get 
\begin{displaymath}
l\bigl(S(r)\bigr) \geq \frac{e^{2r}}{8} \sum_{i=i_0}^\infty \alpha_i^2=C_s^2 \frac{e^{2r}}{8} \sum_{i=i_0}^\infty \frac{1}{i^{2s}} \geq C_s^2 \frac{e^{2r}}{8} \int_{i_0}^\infty \frac{1}{x^{2s}}dx= C_s^2 \frac{e^{2r}}{8(2s-1)i_0^{2s-1}}.
\end{displaymath}

Replacing our explicit value for $i_0$ gives 
\begin{displaymath}
l(S(r)) \geq C e^\frac{r}{s}
\end{displaymath}
for sufficiently large $r$ and some constant $C$ (again depending on $s$). Hence $\underline{\ent} K \geq \frac{1}{s}$.

{\bf Upper bound}\\
For the upper bound in the statement, we apply our main theorem. For this, we have to find an upper bound on the Minkowski dimension of $\ex K=E$. 

Since the Minkowski dimension is invariant under bi-Lipschitz maps, we may replace distances on the unit circle by angular distances. 

$E$ has two accumulation points $\pm x_0$. For $\epsilon>0$, let $N(\epsilon)$ be the number of $\epsilon$-balls needed to cover $E$. We take one such ball around $\pm x_0$ and one further ball for each point in $E$ not covered by these two balls. 

The three points corresponding to the angle $\alpha_i$ are certainly in the $\epsilon$-neighborhood of $\pm x_0$ provided
\begin{displaymath}
3 \sum_{j=i}^\infty \alpha_j \leq \epsilon.
\end{displaymath}

Now we compute that 
\begin{displaymath}
\sum_{j=i}^\infty \alpha_j =C_s \sum_{j=i}^\infty \frac{1}{j^s} \leq C_s \int_{i-1}^\infty \frac{1}{x^s}dx=\frac{C_s}{s-1} \frac{1}{(i-1)^{s-1}}.   
\end{displaymath}

It follows that all $i \geq i_0:=
 \left(\frac{3C_s}{s-1}\right)^\frac{1}{s-1} \epsilon^{\frac{1}{1-s}}+1$ satisfy the inequality above and hence 
\begin{displaymath}
N(\ex K,\epsilon) \leq 6i_0+2 \leq C \epsilon^{-\frac{1}{s-1}}.
\end{displaymath}

It follows that the upper Minkowski dimension is not larger than $\frac{1}{s-1}$. The upper bound of theorem \ref{thm_minkowski_bound} gives 
\begin{displaymath}
\overline{\ent} K \leq \frac{2s-2}{3s-4}.	
\end{displaymath}
\endproof

%---------------------------------------------------------
\section{Centro-projective and centro-affine areas}
\label{sec_centro_proj}

In this section, we will take a closer look at the centro-projective area which was introduced (in a non-intrinsic way) in definition \ref{def_cp_area}. 

\subsection{Basic definitions and properties} 
Geometrically speaking, both centro-affine and centro-projective areas are Riemannian volumes of the boundary $\partial K$. 

We first give intrinsic definitions of the centro-affine metric and area. Let $K$ be a convex body with a distinguished interior point which we may suppose to be the origin $o$ of $V$. The Minkowski functional of $K$ is the unique positive function $F$ that is homogeneous of degree one and whose level set at height $1$ is the boundary $\partial K$. This function is convex and, according to Alexandroff's theorem, has almost everywhere a quadratic approximation.

\begin{definition}
Let $v$ be a tangent vector to $\partial K$ at a smooth point $p$. Then the \emph{centro-affine semi-norm} of $v$ is 
\begin{displaymath}
\Vert v\Vert_a:=\sqrt{Hess_pF(v,v)}.
\end{displaymath}
\end{definition}

The square of the centro-affine semi-norm is a quadratic function on the tangent, hence we may define as usual a volume form, say $\omega_a$ (which vanishes if $\|\cdot\|_a$ is not definite). 

\begin{definition}
The \emph{centro-affine area} of $K$ is
\begin{displaymath}
\cA_a(K):=\int_{\partial K}|\omega_a|.
\end{displaymath}
\end{definition}

It easily follows from the definitions that the centro-affine area is indeed an affine invariant of pointed convex bodies. Moreover, it is finite and vanishes on polytopes. The next proposition relates our definitions with the classical ones, its proof is a straightforward computation.
\begin{proposition}
If the space is equipped with a Euclidean inner pro\-duct, then the centro-affine area is given by
\begin{displaymath}
\cA_a(K)=\int_{\partial K}\frac{\sqrt{k}}{\langle n,p\rangle^{\frac{n-1}{2}}}\ dA,
\end{displaymath}
where $k$ is the Gaussian curvature of $\partial K$ at $p$, $n$ the unit vector normal to $T_p\partial K$ and $dA$ the Euclidean area.
\end{proposition}

In order to introduce the centro-projective area, we will consider a compact convex subset of the (real) $n$-dimensional projective space. Here the word ``convex'' means that each intersection with a projective line is connected. 
 
The definitions of the centro-projective semi-norm and area are merely the same as the centro-affine ones, but one has to replace the Minkowski functional by a projectively invariant function.

\begin{definition}
Let $K\subset \dP^n$ be a convex body and $o \in \inte K$. The \emph{projective gauge function} is
\begin{align*}
G_K\colon \dP^n \setminus \{o\} & \to \dR \cup \{\infty\},\\
x & \mapsto2[q_1,o,x,q_2]
\end{align*}
where $q_1$ and $q_2$ are the two intersections of $\partial K$ with the line going through $o$ and $x$.
\end{definition}

Since the order of $q_1$ and $q_2$ is not fixed, this function is multi-valued (in fact $2$-valued). Identifying $\mathbb{R} \cup \{\infty\}$ with $\mathbb{P}^1$, this function is continuous. 

If $p$ belongs to the boundary of $K$, then the two values of $G_K(p)$ are different, one of them being $2$, the other being $\infty$. Hence there is some neighborhood $U$ of $p$ such that the restriction of $G_K$ to $U$ is the union of two continuous (in fact smooth) functions $G_K^+, G_K^-$ on $U$, where $G_K^+(p)=2$ and $G_K^-(p)=\infty$. 

Let $v$ be a tangent vector to $\partial K$ at a smooth point $p$. Since the restriction of $G_K^+$ to $\partial K \cap U$ is constant, the derivative of $G_K^+$ in the direction of $v$ vanishes. Therefore, the Hessian of the restriction of $G_K^+$ to the tangent line is well-defined. 
 
\begin{definition}
The \emph{centro-projective semi-norm} of $v$ is 
\begin{displaymath}
\Vert v\Vert_p:=\sqrt{Hess_pG_K^+(v,v)}.
\end{displaymath}

Calling $\omega_p$ the induced volume form on $\partial K$, the \emph{centro-projective area} of $K$ is
\begin{displaymath}
\cA_p(K):=\int_{\partial K}|\omega_p|.
\end{displaymath}
\end{definition}

As a consequence of the definition, one has
\begin{proposition}\label{euklid}
In a Euclidean space,
\begin{equation*}
\label{eq_centro_proj_eucl}
\cA_p(K)=\int_{\partial K}\frac{\sqrt{k}}{\langle n,p\rangle^{\frac{n-1}{2}}}\left(\frac{2a}{1+a}\right)^{\frac{n-1}{2}}\ dA.
\end{equation*}
In particular, the intrinsic definition of $\cA_p$ agrees with the definition given in the introduction. 
\end{proposition}

\proof
An easy computation shows that 
\begin{displaymath}
[q_1,o,x,q_2]=\frac{1+a(q_2)}{F(x)+a(q_2)}F(x). 
\end{displaymath}

Then, if $p$ is a smooth point of $\partial K$ and $v\in T_p\partial K$,
\begin{displaymath}
Hess_pG_K(v,v)=\frac{2a(p)}{1+a(p)}Hess_pF(v,v).
\end{displaymath}
\endproof

\subsection{Properties of the centro-projective area} 

Both centro-affine and centro-projective areas vanish on polytopes, hence they are not continuous with respect to the Hausdorff topology on (pointed) bounded convex bodies. Nevertheless, the centro-affine area is upper-semi continuous (see \cite{lut96}). The same holds true for the centro-projective area as shown in the next theorem.
\begin{theorem} \label{thm_properties_ap}
 The centro-projective area is finite, invariant under projective transformations and  upper-semicontinuous.
\end{theorem}
 
\proof
From the above intrinsic definition, it follows that $\cA_p$ is invariant under projective transformations. Also, since the function $a$ on the boundary is bounded and positive and since the centro-affine area is finite, it follows from proposition~\ref{euklid} that the centro-projective area is also finite. It remains to show that it is upper-semicontinuous. Our proof is based on the fact that the centro-affine surface area $\mathcal{A}_a$ is semicontinuous, see E.~Lutwak~\cite{lut96}. 

Let $K$ be a bounded convex body containing the origin in its interior and $(K_i)$ a sequence of convex bodies with the same properties converging to $K$. Set 
\begin{displaymath}
 \tau(p):=\left(\frac{2a(p)}{1+a(p)}\right)^{\frac{n-1}{2}}, \quad p \in \partial K
\end{displaymath}
which is a continuous function on $\partial K$. 

For each $i$, if  $a_i$ is the function corresponding to $K_i$ and $p_i$ is the radial projection of $p$ on $\partial K_i$, 
define $\tau_i \in C(\partial K)$ by 
\begin{displaymath}
 \tau_i(p):=\left(\frac{2a_i(p_i)}{1+a_i(p_i)}\right)^{\frac{n-1}{2}}\text{.}
\end{displaymath}

Since $K_i \to K$, $\tau_i$ converges uniformly to $\tau$. Therefore, for fixed $\epsilon>0$ and all sufficiently large $i$, 
\begin{displaymath}
 \|\tau_i-\tau\|_\infty <\epsilon
\end{displaymath}

Take a triangulation of the sphere and let $\partial K=\cup_{j=1}^m \Delta_j$ (resp. $\partial K_i=\cup_{j=1}^m \Delta_{ij}$) be its radial projection.   

Choosing this triangulation sufficiently thin, there exist $t_1,\ldots,t_m \in \dR_+$ with 
\begin{displaymath}
 |\tau(p)-t_j|<\epsilon
\end{displaymath}
on $\Delta_j$. By the triangle inequality, $|\tau_i(p)-t_j|<2\epsilon$ on $\Delta_{ij}$. 

We define 
\begin{displaymath}
 \mathcal{A}_p(K_i,\Delta_{ij}):=\int_{\Delta_{ij}} \frac{\sqrt{k(x)}}{\bigl\langle n(x),x\bigr\rangle^{\frac{n-1}{2}}} \tau_i d\mathcal{H}^{n-1}(x).
\end{displaymath}
Clearly, $\mathcal{A}_p(K_i)=\sum_{j=1}^m \mathcal{A}_p(K_i,\Delta_{ij})$. In a similar way, we define $\mathcal{A}_p(K,\Delta_{j})$, $\mathcal{A}_a(K_i,\Delta_{ij})$ and $\mathcal{A}_a(K,\Delta_{j})$. 

Fix $p_j$ in the interior of $\Delta_j$ and consider the convex hull $\widehat{\Delta_{i}}$ (resp. $\widehat{\Delta_{ij}}$) of $\Delta_j$ (resp. $\Delta_{ij}$) and $-p_j$. The boundary of $\widehat{\Delta_{ij}}$ is a union of $\Delta_{ij}$ and line segments, hence $\mathcal{A}_a(K_i,\Delta_{ij})= \mathcal{A}_a(\widehat{\Delta_{ij}})$. By the semicontinuity of $\mathcal{A}_a$, we obtain   
\begin{displaymath}
 \limsup_{i \to \infty} \mathcal{A}_a(K_i,\Delta_{ij})= \limsup_{i \to \infty} \mathcal{A}_a(\widehat{\Delta_{ij}})\leq \mathcal{A}_a(\widehat{\Delta_j}) = \mathcal{A}_a(K,\Delta_j). 
\end{displaymath}

It follows that 
\begin{align*}
 \limsup_{i \to \infty} \mathcal{A}_p(K_i) & =\limsup_{i \to \infty} \sum_{j=1}^m \mathcal{A}_p(K_i,\Delta_{ij}) \\
& \leq \limsup_{i \to \infty} \sum_{j=1}^m \mathcal{A}_a(K_i,\Delta_{ij})(t_j+2\epsilon) \\
& \leq \sum_{j=1}^m \mathcal{A}_a(K,\Delta_{j})(t_j+2\epsilon)
\end{align*}

On the other hand, 
\begin{displaymath}
 \mathcal{A}_p(K) = \sum_{j=1}^m \mathcal{A}_p(K,\Delta_j) \geq \sum_{j=1}^m \mathcal{A}_a(K,\Delta_j)(t_j-\epsilon)
\end{displaymath}
from which we deduce that 
\begin{displaymath}
\limsup_{i \to \infty} \mathcal{A}_p(K_i) \leq \mathcal{A}_p(K)+3\epsilon \mathcal{A}_a(K).
\end{displaymath}
\endproof

The centro-affine surface area has the following important properties: 
\begin{enumerate}
\item $\cA_a$ is a valuation on the space of compact convex subsets of $V$ containing $o$ in the interior. This means that whenever $K,L, K \cup L$ are such bodies, then 
\begin{displaymath}
	\cA_a(K \cup L)=\cA_a(K)+\cA_a(L)-\cA_a(K \cap L). 
\end{displaymath}
\item $\cA_a$ is upper semi-continuous with respect to the Hausdorff topology.  
\item $\cA_a$ is invariant under $GL(V)$. 
\end{enumerate}

A recent theorem by M.~Ludwig \& M.~Reitzner \cite{lurei08} states that the vector space of functionals with these three properties is generated by the constant valuation and $\cA_a$. The centro-projective surface area satisfies the last two conditions, but is not a valuation.

%\bibliographystyle{plain}
%\bibliography{biblio}
\end{document}